\documentclass[a4paper,11pt,twoside]{article}
\newcommand{\heuteIst}{Nov 17, 2000 }
\usepackage{amsmath,amsthm,amsfonts,latexsym,amscd,amssymb}

\usepackage{amsthm}

\swapnumbers
\theoremstyle{plain}
\newtheorem{theorem}{Theorem}[section]
\newtheorem{lemma}[theorem]{Lemma}

\newtheorem{proposition}[theorem]{Proposition}

\theoremstyle{definition}
\newtheorem{definition}[theorem]{Definition}
\newtheorem{example}[theorem]{Example}

\theoremstyle{remark}
\newtheorem{remark}[theorem]{Remark}

\usepackage{amsmath,amsfonts,latexsym}

\newcommand{\reals}{\mathbb{R}}

\newcommand{\complexs}{\mathbb{C}}

\newcommand{\naturals}{\mathbb{N}}

\newcommand{\integers}{\mathbb{Z}}

\newcommand{\rationals}{\mathbb{Q}}

\DeclareMathOperator{\id}{id}

\newcommand{\boundedops}{{\mathcal{B}}}

\newcommand{\abs}[1]{\left\lvert#1\right\rvert} 
\newcommand{\norm}[1]{\left\lVert#1\right\rVert}

\newcommand{\tensor}{\otimes}

\newcommand{\disjointunion}{\amalg}
\newcommand{\subgroup}{<}

\DeclareMathOperator{\spec}{Spec}

\DeclareMathOperator{\tr}{tr}
\DeclareMathOperator{\pr}{pr}
\DeclareMathOperator{\coker}{coker}

\DeclareMathOperator{\ind}{ind}

\newcommand{\forget}[1]{}

\newcommand{\innerprod}[1]{\langle #1 \rangle}

{\catcode`@=11\global\let\c@equation=\c@theorem}

\allowdisplaybreaks[2]

\newcommand{\NeumannN}{{\mathcal{N}}}

\begin{document}

\title{The trace on the K-theory of group $C^*$-algebras}
\author{
Thomas Schick\thanks{
e-mail: thomas.schick@math.uni-muenster.de
\protect\\
www:~http://math.uni-muenster.de/u/schickt/
\protect\\
Fax: ++49 251/83 38370
\protect\\
Author's work
funded by Deutscher Akademischer Austauschdienst}\\
FB Mathematik --- Uni M\"unster\\
Einsteinstr. 62 --- 48149 M\"unster, Germany\\
}
\date{Last compiled: \today, last edited  \heuteIst or later}
\maketitle

\begin{abstract}
The canonical trace on the reduced $C^*$-algebra of a discrete group
gives rise to a homomorphism from the K-theory of this $C^*$-algebra
to the real numbers. This paper studies the range of this
homomorphism. For torsion free groups, the Baum-Connes conjecture
together with
Atiyah's $L^2$-index theorem implies that the range consists of the
integers. 

We give a direct and elementary proof that if $G$ acts on a tree and
admits a  homomorphism $\alpha$ to another group $H$
 whose restriction $\alpha|_{G_v}$ to every stabilizer group of a vertex is
 injective, then
 $$\tr_G(K(C_r^*G))\subset \tr_H(K(C_r^*H)).$$
 This follows from a general relative Fredholm
 module technique.

Examples are in particular HNN-extensions of $H$ where the stable
letter acts by conjugation with an element of $H$, or
amalgamated free products $G=H*_U H$ of two copies of the same groups
along a subgroup $U$.

\end{abstract}

\section{Introduction}

Let $G$ be a discrete group. All discrete groups considered in this
paper are assumed to be
countable. The trace $\tr_G\colon \complexs G\to\complexs\colon 
\sum_{g\in G} \lambda_g g\mapsto \lambda_1$ (where $1$ is the neutral
element of $G$) extends to a trace on the reduced $C^*$-algebra of
$G$ and therefore gives rise to a homomorphism
\begin{equation*}
  \tr_G\colon  K_0(C^*_{r}G)\to\reals.
\end{equation*}
If $G$ is torsion free, we have the commutative diagram
\begin{equation*}
  \begin{CD}
    K_0(BG) @>A>> K_0(C^*_{r}G)\\
    @VV{\ind_G}V  @VV{\tr_G}V\\
    \integers @>>> \reals ,
  \end{CD}
\end{equation*}
where $A$ is the Baum-Connes assembly map.
The Baum-Connes conjecture says that $A$ is an
isomorphism. We denote
by $\ind_G$ Atiyah's $L^2$-index, which coincides with the ordinary
index \cite{Atiyah(1976)} and therefore takes
values in the integers. Surjectivity of $A$ of course implies that
$\tr_G$ is also
integer valued. We will denote this consequence of the Baum-Connes
conjecture as the \emph{trace conjecture}. It implies by a standard
argument that there
are no nontrivial
projections in $C^*_{r}G$. 

The trace conjecture was  verified directly for free
groups using a special Fredholm module, which can be
assigned to free groups, cf.~e.g.~\cite{Effros(1989)}. Based on the ideas of this proof we
get the following result.

\begin{theorem}\label{Ktraceprop}
    Let $H, G$ be discrete countable groups and assume
    \begin{equation*}
\tr_H(K(C^*_rH))\subset A\subset\reals .
\end{equation*}
 Let $\Omega$ and $\Delta$ be
  sets with commuting $G$-action from the left and $H$-action from the
  right such that $\Omega$ and $\Delta$ are free $H$-sets. Let
  $\Omega=\Omega'\cup X$ and assume $\Delta$ and $\Omega'$ are free
  $G$-sets  and
  $X$ consists of $1\le r<\infty$ $H$-orbits. Assume
  there is a bijective right $H$-map $\phi\colon \Delta\to\Omega$. Suppose
  that
  for every $g\in G$ the set
  \begin{equation*}
    R_g:=\{ x\in \Delta\mid \phi(gx)\ne g\phi(x)\}
  \end{equation*}
  is contained in the union of finitely many $H$-orbits of $\Delta$.

  Then 
  \begin{equation*}
    \tr_G(K(C^*_rG))\subset \frac{1}{r} A.
  \end{equation*}
\end{theorem}

It remains now of course to give interesting examples of Theorem
\ref{Ktraceprop}. For this, we use results of Dicks-Schick
\cite{Dicks-Schick(2000)} and the following definition:

\begin{definition}
 Let $G$ be a group which acts on a tree. We say $G$ is \emph{subdued} by a
  group $H$ if
  there is a homomorphism $\alpha\colon G\to H$ whose restriction
  $\alpha|_{G_v}$ to the stabilizer group $G_v:=\{g\in G\mid gv=v\}$ is
  injective
  for each vertex $v$ of the tree.
\end{definition}

Remember that one can translate between groups acting on trees and
fundamental groups of graphs in such a way that conjugacy classes of vertex
stabilizers correspond to vertex groups.

\begin{example}\label{ex:subdue}
  \begin{enumerate}
  \item\label{item:HNN} If $G=H*_U H$, then $G$ is the fundamental group
    of a graph 
    consisting of two vertices joined by an edge. Hence there is an
    action of $G$ on a tree such that the stabilizer group of every
    vertex is conjugate to one of the two copies of $H$. Using the
    obvious projection $G\to H$ which is the identity on both factors
    we see that $G$ is subdued by $H$.
  \item \label{item:amalg} Assume $U\subgroup H$ and $g\in H$. Let
    $\phi\colon U\to U^g$ be 
    given by conjugation with $g$. Let $G$ be the HNN-extension of $H$
    along $\phi$. Then $G$ is the fundamental group of a graph of
    groups with one vertex (the vertex group being $H$) and therefore
    acts on a tree such that each vertex stabilizer is conjugate to
    $H$. We define a homomorphism $\alpha\colon G\to H$ such that the
    restriction of $\alpha$ to $H$ is the identity and  maps the
    stable letter $t$ to $g$. This implies that $G$ is subdued by $H$.
\end{enumerate}
\end{example}

\begin{theorem}\label{heredtrace}
Assume  $G$ is the fundamental group of a graph of groups which
  is subdued by another group $H$. Then
  \begin{equation*}
    \tr_{G}\left( K_0(C^*_{r}[G])\right) \subset 
  \tr_{H} \left(K_0(C^*_{r}[H])\right) .\end{equation*}
  In particular, this applies to the situations of Example \ref{ex:subdue}.
\end{theorem}
If $\tr_H(K_0(C^*_r[H]))\subset
\integers$ then this
implies the trace conjecture for $G$. Baum and Connes
conjecture \cite[p.32]{Baum-Connes(1982)} that for
groups with torsion the range should be contained in $\rationals$.
Our results support also this assertion (Baum and Connes make in fact a more
precise and
stronger conjecture, which however was disproved by Roy
\cite{Roy(1998)}).

\begin{remark}\label{remark:otherway}
  In some cases, it is possible to derive  the conclusions of Theorem
  \ref{heredtrace} from elaborate K-theory calculations. One can use
  the exact  sequence for
the fundamental group of a graph of groups \cite[Theorem
  18]{Pimsner(1986)}. In some cases elementary
  properties of the trace then imply Theorem \ref{heredtrace}.

  One
  might hope to give a general treatment of the range of the traces in
  these cases as is done for certain HNN-extensions in
  \cite{Pimsner(1985)} and \cite{Exel(1987)}. However, even in the
  case of HNN-extensions, in
  general those results are difficult
  to interpret and it is not clear that \cite{Pimsner(1985)} or
  \cite{Exel(1987)} 
  implies Example \ref{ex:subdue}~\ref{item:HNN}.

  Moreover, observe
  that Pimsner uses deep
  KK-theoretic methods to derive the
 exact sequence for a graph of groups. 
In contrast, our derivation is elementary.

\end{remark}

\begin{remark}
  To apply Theorem \ref{heredtrace}, essentially we have to know the
  trace
  conjecture for $H$. The obvious sufficient condition is that $H$
  fulfills the Baum-Connes conjecture, but it is also enough that $H$
  is a subgroup of such a group.

  If $H$ satisfies the Baum-Connes conjecture with coefficients, then
  the same is true for each group $G(v)$
  (since they are subgroups of $H$),
  hence by \cite{Oyono-Oyono(1998)} $G$ fulfills the Baum-Connes
conjecture, and the statement of Theorem \ref{heredtrace} follows also
immediately from this fact.

However, these arguments do not apply to the Baum-Connes conjecture without
coefficients. Lafforgue proves that cocompact discrete subgroups
e.g.~of $Sl_3(\reals)$ or $Sp(n,1)$ satisfy the Baum-Connes conjecture
without coefficients \cite{Lafforgue(1998)}. However, it is unknown
whether the Baum-Connes 
conjecture with coefficients is true for these groups. Therefore, the
  consequences of Theorem
  \ref{heredtrace} are not included
  in the knowledge about the Baum-Connes conjecture for a non-cocompact
  subgroup $H$ of $Sl_3(\reals)$ or $Sp(n,1)$ contained in a cocompact
  torsion-free subgroup.

  In particular, for such an $H$ every $H*_U H$ fulfills the trace
  conjecture, but it is not
  clear whether it fulfills the Baum-Connes conjecture.
\end{remark}

The method described in \cite{Effros(1989)} was used by Linnell
\cite{Linnel(1993),Linnell(1998)} to prove the Atiyah conjecture about
the
integrality of $L^2$-Betti numbers for free groups, and starting with
this for a lot of other groups. We investigate the Atiyah conjecture and
obtain  generalizations of Linnell's results in
\cite{Schick(1999),Dicks-Schick(2000)}.

\section{The trace conjecture for the K-theory of group
  $C^*$-algebras}
\label{sec:trace}

In this section we prove Theorem \ref{Ktraceprop} and Theorem
\ref{heredtrace}.

We first show how Theorem \ref{heredtrace} follows from
Theorem \ref{Ktraceprop}. The method for this is developed by
Dicks and Schick in
\cite{Dicks-Schick(2000)}. For
the convenience of the reader we repeat
the easy proof of the special case we are concerned with here.

\begin{proof}[Proof of Theorem \ref{heredtrace}]
Assume $G$ acts on the tree $T$ with set of vertices $V$ and set of
edges $E$.    We choose an arbitrary $v_0 \in V$.

   Let $\{\ast\}$ be a trivial $G$-set.  Let
   $\tilde\phi\colon V \to E\,\vee\,\{\ast\}$ denote the map which assigns to
   each $v \in V$ the last edge in the $T$-geodesic from $v_0$ to $v$,
   where this is taken to be $\ast$ if $v = v_0$.  By
   Julg-Valette~\cite{Julg-Valette(1984)}, $\tilde\phi$ is bijective, and,
   for all $v \in V$, $g \in G$, we have $\tilde\phi(gv) =
   g\tilde\phi(v)$ if and
   only if $v$ is not in the $T$-geodesic from $v_0$ to $g^{-1}v_0$.

   Define now $\Delta:=V\times H$, $\Omega':=E\times H$ and
   $\Omega:=(E\disjointunion \{\ast\})\times H$. We define the $G$ and
   the $H$-action on $\Delta$ and $\Omega'$ be setting
   \begin{equation*}
     g(x,u)h:= (gx,\alpha(g)uh)\qquad\forall g\in G,\;x\in T\; u,h\in H.
   \end{equation*}
   Since the restriction of $\alpha$ to each stabilizer group is
   injective, this is a free $G$- and of course also a free $H$-action,
   and they commute. Extend the action to $\Omega$ by
   $g(\ast,u)h=(\ast,uh)$ for $g\in G$ and $u,h\in H$.

   We define $\phi\colon \Delta=V\times
   H\to\Omega=(E\disjointunion\{\ast\})\times H$ by
   \begin{equation*}
     \phi(v,u)= (\tilde\phi(v),u).
   \end{equation*}
   Since $\tilde\phi$ is bijective, the same is true for
   $\phi$. Moreover, for fixed $g\in G$ we have
   $\phi(g(v,u))=(\phi(gv),\alpha(g)u) = (g\phi(v),\alpha(g)u)=
   g\phi(v,u)$ exactly if $\phi(gv)= g\phi(v)$. Therefore, the
   set $R_g$ of Theorem \ref{Ktraceprop}
   is contained in the union of the finitely many $H$-orbits of
   $\Delta$ determined by the $T$-geodesic from $v_0$ to
   $g^{-1}v_0$. Theorem \ref{heredtrace} now follows from Theorem
   \ref{Ktraceprop}. 
\end{proof}

We are now going to prove Theorem \ref{Ktraceprop}. 
We use the language of Hilbert $A$-modules (for a $C^*$-algebra $A$),
compare e.g.~\cite[Section 13]{Blackadar(1986)}.
\begin{definition}\label{defofE}
  Let $T$ be a free set of generators of the free $H$-set $\Delta$.
  Then we can form the Hilbert $C^*_rH$-module
  $E:=l^2(T)\tensor_\complexs C^*_r H$. Of course, $E$ is nothing but
  the Hilbert sum of copies of $C^*_rH$ indexed by $T$ (in the sense
  of Hilbert $C^*_rH$-modules). Moreover, $E\tensor_{C^*_r H} l^2
  H\cong l^2(\Delta)$, and $E$ is in a natural way a subset of
  $l^2(\Delta)$ because $C^*_rH\subset l^2 H$.  The module $E$ as a
  subset of $l^2(\Delta)$ does not depend on the choice of the basis
  $T$.
\end{definition}

Let $\boundedops(E)$ be the set of bounded adjointable Hilbert
$C^*_rH$-module homomorphisms. The map 
\begin{equation*}
  \boundedops(E)\to \boundedops(l^2(\Delta))\colon  A\mapsto A_\Delta:=A\tensor 1
\end{equation*}
is an injective algebra homomorphism by
\cite[p.~111]{Blackadar(1986)}, therefore an isometric injection of
$C^*$-algebras \cite[4.1.9]{Kadison-Ringrose(1983)}. Observe that the
image of $\boundedops(E)$ commutes with the right action of $H$ on
$l^2(\Delta)$ and therefore is contained in the corresponding von
Neumann algebra with its canonical trace. We will prove in Lemma
\ref{inclusion} that
$\complexs G\subset \boundedops(E)$. It follows that the closure of
$\complexs G$ in $\boundedops(E)$ and $\boundedops(l^2(\Delta))$
coincides. Since $\Delta$ is a free $G$-set, this closure is
isomorphic to $C^*_rG$.

\begin{lemma}\label{inclusion}
  We have a natural inclusion $\complexs G\subset\boundedops(E)$ which
  extends the action of $G$ on $\complexs\Delta\subset E$, in
  particular it is compatible with the injection $\complexs
  G\subset\boundedops(l^2\Delta)$.
\end{lemma}
\begin{proof}
  First observe that $E$ is a closure of the algebraic tensor product
  of $l^2(T)$ 
and $C^*_rH$. Moreover, $\complexs[T]$
  is dense in $l^2(T)$ and
  $\complexs H$ is dense in $C^*_rH$. Therefore
  $\complexs\Delta=\complexs[T]\tensor_\complexs\complexs H$ is a
  dense subset of $E$.

  Now fix $g\in G$. Since $\Delta=\bigcup_{t\in T} tH$, where the union
  is disjoint, for each $t\in T$ we get unique elements $t_{g,t}\in
  T$, $h_{g,t}\in
  H$ such that $gt=t_{g,t}h_{g,t}$. Since the actions of $G$ and $H$
  commute, $gt=gt'h$ implies $t=t'$. The map $\alpha_g\colon T\to
  T\colon t\mapsto
  t_{g,t}$ therefore is a bijection. Pick 
  \begin{equation*}
x=\sum_{t\in T}tv_t, \;x'=\sum_{t\in T}t v'_t\in
  \complexs[T]\tensor_{\complexs} \complexs H, \quad\text{with }v_t,v'_t\in
  \complexs H.
\end{equation*}
By linearity we get 
\begin{equation*}
gx=\sum_{t\in T}
  t_{g,t}h_{g,t}v_t\quad\text{and}\quad gx'=\sum_{t\in T}
  t_{g,t}h_{g,t}v'_t.
\end{equation*}
As is the convention in the theory of Hilbert $A$-modules,  all of our
    inner products are
    linear in the second variable. 
Taking now the $C^*_rH$-valued inner product we 
  get (adjoint and products of elements in $\complexs H\subset C^*_r
  H$ are taken in the sense of the $C^*$-algebra)
  \begin{equation*}
    \begin{split}
      \innerprod{x,x'}_{C^*_rH} & = \sum_{t\in T}v_t^*v_t' = \sum_{t\in T}
      (h_{g,t}v_t)^* (h_{g,t}v_t')\\
      &= \sum_{t\in T} (h_{g,\alpha_g^{-1}(t)}v_{\alpha_g^{-1}(t)})^*
      (h_{g,\alpha_g^{-1}(t)}v'_{\alpha_g^{-1}(t)}) = \innerprod{gx,gx'}_{C^*_rH}
    \end{split}
    \end{equation*}
  (here we used that $H$ acts unitarily,
    i.e.~$h_{g,t}^*=h_{g,t}^{-1}$). Hence $G$ acts
    $C^*_rH$-isometrically and this action extends to $E$. By
    linearity we get an $*$-algebra homomorphism  $\complexs
    G\to\boundedops(E)\to \boundedops(l^2\Delta)$. The composition
    is injective, therefore the same is true for the first map.
\end{proof}

The above reasoning implies in the same way:
\begin{lemma}
  For $n\in\naturals$ we have canonical injections of $*$-algebras
  \begin{equation*}
    M_n(\complexs G)\subset M_n(C^*_r G)\subset
    \boundedops(E^n)\subset \boundedops(l^2\Delta^n)^H,
  \end{equation*}
  where $\boundedops(l^2\Delta^n)^H$ denotes the operators which commute
  with the right action of $H$.
\end{lemma}

We have to compute the $G$-trace of operators in $M_n(C^*_r G)$, and
we want to express this in terms of the $H$-trace of a suitable other
operator. For this end, we repeat the following definition:
\begin{definition}\label{deftrH}
  An operator $A\in\boundedops(E)$ is of $H$-trace class if its image
  $A_{l^2\Delta}\in\boundedops(l^2(\Delta))^H$ is of $H$-trace  class in the
  sense of the von Neumann algebra  (compare
  e.g.~\cite[2.1]{Schick(1998c)}), i.e.~if (with $\abs{A}=\sqrt{A^*A}$)
  \begin{equation*}
    \sum_{t\in T} \innerprod{t,\abs{A_{l^2\Delta}}t}_{l^2(\Delta)}< \infty.
  \end{equation*}
Then also $ \sum_{t\in T} \innerprod{t,A_{l^2\Delta}t}_{l^2(\Delta)}$ converges
and we set
  \begin{equation*}
    \tr_H(A):= \sum_{t\in T} \innerprod{t,A_{l^2\Delta}t}_{l^2(\Delta)}.
  \end{equation*}
  For $A=(A_{ij})\in\boundedops(E^n)=M_n(\boundedops(E))$
 we set
  \begin{equation*}
    \tr_H(A)= \sum_{i=1}^n \tr_H(A_{ii}),
  \end{equation*}
    if $\abs{A_{l^2\Delta}}$ of
  $H$-trace class (with the obvious definition for this).
\end{definition}

\begin{lemma}\label{computetrH}
  Let $A,B,C\in\boundedops(E^n)$ and $A$ be
  of $H$-trace class.
  The trace class operators form an ideal inside $\boundedops(E^n)$ and
  we have $\tr_H(AB)=\tr_H(BA)$. Moreover
  \begin{align*}
    \abs{\tr_H(A+B)} \le & \tr_H(\abs{A}) + \tr_H(\abs{B})\\
    \tr_H({\abs{CAB}})\le & \norm{C}\cdot\norm{B}\cdot\tr_H(\abs{A}).
  \end{align*}
\end{lemma}
\begin{proof}
  These are standard properties of the (von Neumann) trace, compare
  \cite[2.3]{Schick(1998c)}, \cite[Th\'eor\`eme 8 and Corollaire 2
  on p.~106]{Dixmier(1969)}.
\end{proof}

Set $S:=\phi(T)$. This is an $H$-basis for the free $H$-set
$\Omega$. Similarly to
$E$ we can build $F:=l^2(S)\tensor_\complexs C^*_rH\subset
l^2(\Omega)$.
If $S':=S\cap\Omega'$ and $S'':=S\cap X$ (i.e.~$S'$ is an $H$-basis for
$\Omega'$ and $S''$ is an $H$-basis for $X$) then with
$F':=l^2(S')\tensor_\complexs C^*_rH$ and
$F'':=l^2(S'')\tensor_\complexs C^*_rH$ we get a direct sum
decomposition of Hilbert $C^*_rH$-modules
\begin{equation*}
 F= F'\oplus F''.
\end{equation*}
As in the case of $E$ we get an canonical inclusion
\begin{equation*}
  \complexs G\subset C^*_rG\subset \boundedops(F')\subset
  \boundedops(F)\subset \boundedops( l^2\Omega)
\end{equation*}
(we extend the action of $C^*_rG$ to all of $F$ by setting it zero on
$F''$). This composition is a non-unital $*$-algebra homomorphism.

Corresponding statements hold for matrices.

Denote the image of $A\in M_n(C^*_rG)$ in $\boundedops(E^n)$ with
$A_\Delta$ and in $\boundedops(F^n)$ with $A_\Omega$. We therefore have
$A_\Omega=A_{\Omega'}\oplus 0$ with the obvious
notation. 

The bijection $\phi\colon \Delta\to\Omega$ induces a unitary
map of Hilbert spaces $\phi\colon l^2(\Delta)^n\to l^2(\Omega)^n$. Since
$\phi\colon \Delta\to \Omega$
is $H$-equivariant, the same is true for the unitary map. Moreover, we get
a Hilbert $C^*_rH$-module unitary map $\phi\colon E^n\to F^n$. One
key observations is now (this is an extension of the corresponding
observation in the classical proof for the free group):

 \begin{lemma}\label{GtoHlemma}
    Suppose $A\in M_n(C^*_r G)\subset \boundedops(E^n)$ is such that
    the Hilbert
    $C^*_r H$-module morphism $A_\Delta - \phi^*A_\Omega \phi\colon E^n\to E^n$ is of
    $H$-trace class. Then
    \begin{equation*}
      \tr_{G}(A)= \frac{1}{r} \tr_{H}( A_\Delta - \phi^*A_\Omega \phi) .
    \end{equation*}
  \end{lemma}
\begin{proof}
  Observe that $\phi$ is diagonal and traces are the sum
  over the diagonal entries. Therefore we may assume that $n=1$.
  Since $\Delta$ is a free $G$-module, for every $x\in \Delta$ (which
  we identify with the element of $l^2(\Delta)$ which is $1$ at
  $x$ and zero everywhere else) we have
  \begin{equation*}
    \innerprod{x,A_\Delta x}_{l^2(\Delta)} = \tr_{G}(A) 
  \end{equation*}
  (simply identify $Gx$ with $G$ and the left and right hand side
  become identical). Similarly, since $A_\Omega= A_{\Omega'}\oplus 0$
  on $l^2(\Omega')^n\oplus l^2(X)^n$
  \begin{equation*}
    \innerprod{\phi x,A_\Omega \phi x}_{l^2(\Omega)}=\begin{cases}
    \tr_{G}(A); & \text{if }\phi(x)\in\Omega'\\
      0; & \text{if } \phi(x)\in X .\end{cases}
  \end{equation*}
Moreover 
  \begin{equation*}
    \begin{split}
      \tr_{H}( A_\Delta - \phi^*A_\Omega \phi) & = \sum_{t\in T}
      \innerprod{t,A_\Delta t}_{l^2(\Delta)} - \innerprod{t,\phi^*
      A_\Omega \phi(t)}_{l^2(\Delta)} \\
     & = \sum_{\phi(t)\in X\cap S''}
      \tr_{G}(A) = r\tr_{G}(A)
\end{split}
\end{equation*}
since $\abs{X\cap S''}$ is the number of $H$-orbits in $X$, i.e.~$r$.
All other summands cancel each other out.
\end{proof}

Because $K(C^*_rG)$ is generated by projections  $P\in
M_n(C^*_rG)\subset M_n(\NeumannN G)$ and the trace we have to compute
is exactly $\tr_G(P)$, we are
tempted to apply Lemma \ref{GtoHlemma} to such a $P$. A problem
is that it is hard to check whether the trace class condition is
fulfilled in general.

To circumvent these difficulties recall the following fact (compare
e.g.~\cite[III.3, Proposition 3]{Connes(1994)}):
\begin{proposition}\label{holomclosed}
  Let $B$ be a $C^*$-algebra and $U\subset B$ a dense $*$-subalgebra that
  is closed under holomorphic functional calculus. Then the inclusion
  induces an isomorphism
  \begin{equation*}
    K(U)\cong K(B).
  \end{equation*}
  In particular, if $B=C^*_rG$ and $U$ is closed under holomorphic
  functional calculus and contains $\complexs G$, then the ranges of
  the canonical trace applied to $K(C^*_rG)$ and $K(U)$ coincide.
\end{proposition}

As algebra $U$ we will use the closure under holomorphic functional
 calculus of $\complexs G\subset\boundedops(E)$. This of course
 fulfills the conditions of Proposition \ref{holomclosed}. It remains
 to check:

 \begin{lemma}
   \label{Utrclass}
   Let $x\in M_n(U)\subset M_n(C^*_r G)\subset \boundedops(E^n)$, where
   $U$ is the
   closure under holomorphic functional calculus of $\complexs G$ in
   $C^*_rG$. Then $x-\phi^*x\phi\colon E\to E$ is of $H$-trace class and
   $C^*_rH$-compact. 
 \end{lemma}
\begin{proof}
  Start with $g\in
 G\subset\complexs G$. Since $R_g$ as defined in Theorem
 \ref{Ktraceprop} is contained
in finitely many $H$-orbits, $gt=\phi^*g\phi t$ for all but a finitely 
many
 $t\in T$. In particular $g-\phi^*g\phi$ is zero
 outside the $C^*_rH$-submodule of $E$ spanned by this finite number
 of elements of $T$, and is nonzero only on the complement, which is
 isomorphic to $(C^*_rH)^N$ for some $N\in\naturals$. Since
 $\id\colon C^*_rH^N\to C^*_rH^N$ is of finite rank in the sense of Hilbert
 $C^*_rH$-module morphisms, the same is true for
 $g-\phi^*g\phi$. Finite rank operators form a subspace, therefore the
 same is true if we replace $g$ by $v\in\complexs G\subset
 \boundedops(E)$. Passage to finite matrices does preserve the finite
 rank  property. Finite rank implies $H$-trace class and
 $C^*_rH$-compactness. In particular, all operators in
 $M_n(\complexs G)$ give
 rise to $C^*_r H$-compact operators, which also are of $H$-trace
 class. The map $x\mapsto x-\phi^*x\phi$ is norm continuous and the
 compact operators form a closed ideal, therefore $x-\phi^*x\phi$ is
 compact even for arbitrary $x\in M_n(C^*_r G)$.

 Assume $A\in M_n(C^*_r H)$ and $0\ne\xi\notin\spec(A)$. Since the
 homomorphism $A\mapsto A_\Delta$ is unital, we get
 \begin{equation*}
\left((\xi-A)^{-1}\right)_\Delta= (\xi-A_\Delta)^{-1}.
\end{equation*}
Similarly
 $\left((\xi-A)^{-1}\right)_{\Omega'} =
 (\xi-A_{\Omega'})^{-1}$. Consequently
 \begin{equation*}
\left( (\xi-A)^{-1}\right)_\Omega
 = \left( (\xi-A)^{-1}\right)_{\Omega'}\oplus 0 = (\xi-A_\Omega)^{-1}
 - \xi^{-1}P
\end{equation*}
where $P\colon F\to F$ is the projection onto $F''$. Note that
$(\xi-A_\Omega)^{-1}$ acts by multiplication with $\xi^{-1}$ on
$F''$. Here we need the assumption $\xi\ne 0$. Since $T''$ is finite,
$P$ is of finite rank as Hilbert
$C^*_rH$-module morphism.

Suppose now $f$ is a function that is holomorphic in a neighborhood
of $\spec(A)$. Let $\Gamma$ be a loop around $\spec(A)$ and choose
$\Gamma$ so that it does not meet $0\in\complexs$. Then
\begin{equation*}
  \begin{split}
    f(A)_\Delta = & \left(\int_\Gamma
      f(\xi)(\xi-A)^{-1}\;d\xi\right)_\Delta
    =  \int_\Gamma f(\xi)(\xi-A_\Delta)^{-1}\;d\xi\\
    f(A)_\Omega = & \int_\Gamma
      f(\xi)\left((\xi-A_\Omega)^{-1}-\xi^{-1}P\right)\;d\xi. 
\end{split}
\end{equation*}
For $u,v\in\boundedops(E)$ and $\xi\notin\spec(u)\cup\spec(v)$ we have
\begin{equation*}
  \begin{split}
    (\xi-u)^{-1}-(\xi-v)^{-1}
    = &(\xi-u)^{-1}(\xi-v-(\xi-u))(\xi-v)^{-1}\\
    = &
    (\xi-u)^{-1}(u-v)(\xi-v)^{-1}.
\end{split}
\end{equation*}
Therefore
  \begin{multline*}
    f(A)_\Delta-\phi^*f(A)_\Omega\phi =\phi^*P\phi\underbrace{\int_{\Gamma}
     f(\xi)\xi^{-1}\;d\xi}_{\in\complexs}\\
    + \int_\Gamma
    \underbrace{f(\xi)(\xi-A_\Delta)^{-1}}_{=:f_\Delta(\xi)}
     \underbrace{(A_\Delta-\phi^*A_\Omega\phi)}_{=:A_0}
     \underbrace{(\xi-\phi^*A_\Omega\phi)^{-1}}_{=:f_\Omega(\xi)}\;d\xi.
  \end{multline*} 

As a consequence of Lemma \ref{computetrH} we have
\begin{equation*}
  \begin{split}
    \tr_H \abs{\int_\Gamma f_\Delta(\xi) A_0 f_\Omega(\xi)\;d\xi} \le &
    \int_\Gamma \tr_H(\abs{ f_\Delta(\xi) A_0 f_\Omega(\xi)})\; d\xi \\
    \le &
    \int \norm{f_\Delta(\xi)}\cdot \norm{f_\Omega(\xi)} \cdot
    \tr_H(\abs{A_0}).
\end{split}
\end{equation*}
Since the operator valued functions $f_\Delta$ and $f_\Omega$ are
norm-continuous, $f(A)_\Delta-\phi^*f(A)_\Omega)\phi$
 is of $H$-trace class if $A_0$ is of $H$-trace class, in particular
 if $A\in M^n(\complexs G)$. This concludes the proof.
\end{proof}

We determine now the range of the trace on the dense and
holomorphically closed subalgebra $U$ of $C^*_r G$. Since the trace
class condition is fulfilled, it only remains to calculate
$\tr_H(P_\Delta-\phi^*P_\Omega\phi)$ for a projection over $U$, and
Theorem \ref{Ktraceprop} follows from Lemma \ref{GtoHlemma}.

\begin{lemma}
  \label{traceindex} 
  Let $E$ be the Hilbert $C^*_rH$-module introduced in Definition
   \ref{defofE} and
   $P,Q\in \boundedops(E)$ be projections such that $P-Q$ is of
  $H$-trace class and compact in the sense of Hilbert $C^*_rH$-module
   morphisms. Then 
   \begin{equation*}
w:=\left( (PE\oplus QE), 1,\left(\begin{smallmatrix}0 & PQ\\ QP &
   0\end{smallmatrix}\right)\right)
\end{equation*}
 is a Kasparov triple (in the sense of \cite[17.1.1]{Blackadar(1986)})
   representing an
   element in $KK(\complexs,C^*_rH)\cong K_0(C^*_rH)$ and 
   \begin{equation*}
     \tr_H(P-Q) = \ind_H(w) = \tr_H([w]),
   \end{equation*}
   where $\tr_H(P-Q)$ is to be understood in the sense of Definition
   \ref{deftrH}, whereas $\tr_H([w])$ is the canonical trace defined on
   $K_0(C^*_rH)$. 
\end{lemma}
\begin{proof}
  Using Lemma \ref{computetrH} and the fact that $P^2=P$, $Q^2=Q$, $(1-P)P=0$,
  and $\tr_H(XY)=\tr_H(YX)$,  we conclude
  \begin{equation*}
    \begin{split}
      \tr_H(P-Q) & =  \tr_H(P^2(P-Q)) + \tr_H((1-P)(P-Q^2))\\
     &=  \tr_H(P(P-Q)P)
    -      \tr_H(Q(1-P)Q)\\
      & =  \tr_H(P- PQP) - \tr_H(Q-QPQ).
    \end{split}
  \end{equation*}
Let $\alpha_P\colon PE\to PE$ be the orthogonal projection with image
$PE\cap\ker(QP)$, and $\alpha_Q\colon QE\to QE$ the orthogonal
projection with image $QE\cap \ker(PQ)$. Observe that
$\alpha_P=(P-PQP)\alpha_P$. Therefore $\alpha_P$ is of $H$-trace
class, since the same is true for $P-PQP$. In the same way we see that
$\alpha_Q$ is of $H$-trace class. Set
\begin{equation*}
  \begin{split}
    T_0 &:= \id_{PE}-PQP -\alpha_P\colon PE\to PE\\
    T_1 &:= \id_{QE}-QPQ -\alpha_Q\colon QE\to QE.
  \end{split}
\end{equation*}
Then
\begin{equation}
  \label{eq:single_out_kernel}
  \begin{split}
    \tr_H(P-Q) &= \tr_H(P-PQP)-\tr_H(Q-QPQ)\\
      &= \tr_H(T_0)-\tr_H(T_1) +\tr_H(\alpha_P)-\tr_H(\alpha_Q).
\end{split}
\end{equation}
Now $QP\colon PE\to QE$ is a bounded operator with adjoint $PQ\colon
QE\to PE$, and 
\begin{equation}
  \label{eq:ker_formula}
  \begin{split}
    \tr_H(\alpha_P) &= \dim_H(\ker(QP\colon PE\to
    QE))\\
    \tr_H(\alpha_Q) &= \dim_H(\ker(PQ\colon QE\to PE))\\
    &=
    \dim_H(\coker(QP\colon PE\to QE)).
\end{split}
\end{equation}
For a complemented submodule $X$, one defines
$\dim_H(X):=\tr_H(\pr_X)$, where $\pr_X$ is the orthogonal projection
onto $X$.
Since $QP\alpha_P=0$ and $PQ\alpha_Q=0$, and the latter implies
$0=(PQ\alpha_Q)^*=\alpha_Q QP$, we have
\begin{multline*}
  QPT_0= QP(\id_{PE}-PQP-\alpha_P)\\
     =\id_{QE} QP-QPQ^2P-\alpha_Q QP = T_1QP.
\end{multline*}
Moreover, $\ker(QP\colon PE\to QE)\subset \ker(T_0)$, since $QP
(Px)=0$ implies $\alpha_P(Px)=Px$, and in the same way we conclude
$\ker(PQ)=\ker((PQ)^*)\subset \ker(T_1)=\ker(T_1^*)$. It follows that
 $QP$ ``conjugates'' $T_0$ and $T_1$, and by \cite[Proposition
 2.6]{Schick(1998c)}
 (which goes back to a corresponding result in
 \cite[p.~67]{Atiyah(1976)}) that 
 \begin{equation*}
   \tr_H(T_0)=\tr_H(T_1).
 \end{equation*}
 Using Equation \eqref{eq:single_out_kernel} and Equation
 \eqref{eq:ker_formula} we arrive at
  \begin{equation*}
    \tr_H(P-Q) = \ind_H(QP\colon  PE\to QE)
  \end{equation*}
  with the obvious definition of $\ind_H$. This is exactly
  the $H$-index in the graded sense of the operator $F:=\left(\begin{smallmatrix}0 & PQ\\ QP &
   0\end{smallmatrix}\right)\colon  PE\oplus QE\to PE\oplus QE$ (where $PE$ is the positive and
  $QE$ the negative part of the graded Hilbert $C^*_rH$-module
  $PE\oplus QE$).

  The only thing it remains to check is whether $w$ fulfills all
  the axioms of Kasparov triples. Since the action of $\complexs$ is
  unital and the operator is self adjoint, this amounts to check that
  $1-F^*F$ and $1-FF^*$ are compact in the sense of Hilbert
  $C^*_rH$-module morphisms. Now $F^*F=F^2=FF^*=
\left(  \begin{smallmatrix}PQP & 0\\ 0 & QPQ 
   \end{smallmatrix}\right)$. Since $P-Q$ is compact, the same is true
  for $P(P-Q)P=P-PQP\colon E\to E$. Then also the composition with the
  inclusion of $PE$ into $E$ and the projection $P\colon E\to PE$ is
  compact. This operator coincides with $1-PQP\colon PE\to PE$. Similarly
  $1-QPQ\colon QE\to QE$ is compact. This concludes the proof.
\end{proof}

To finish the proof of Theorem \ref{Ktraceprop} observe that by
Proposition \ref{holomclosed} it suffices to compute $\tr_G(P)$ if
$P\in M_n(U)$ is a projection, where $U$ is the holomorphic closure of
$\complexs G\subset C^*_rG$. Since $A\to A_\Delta$ and $A\to A_\Omega$
are $*$-algebra homomorphisms, $P_\Delta$ and $P_\Omega$ are
projections. Now Lemma \ref{Utrclass} implies that we
can apply Lemma \ref{traceindex} to $P_\Delta-\phi^*P_\Omega\phi$. By
assumption $\tr_H(K_0(C^*_rH))\subset A$, therefore
$\tr_H(P_\Delta - \phi^*P_\Omega\phi)\in A$. By Lemma \ref{GtoHlemma}
then
$\tr_G(P)\in\frac{1}{r} A$, and this concludes the proof of
Theorem \ref{Ktraceprop}.

\begin{remark}
  We use the language of Hilbert modules and Kasparov triples only for
  convenience. Observe that we don't use much more than the
  definition: the single theorem we use is that our Kasparov triples
  indeed give rise to K-theory elements, and this is not very
  deep. By \cite[17.5.5]{Blackadar(1986)} $KK(\complexs,C^*_rG)$ and
  $K_0(C_r^*G)$ are isomorphic, but to construct the map much less is
  needed. (Essentially we only have to perturb $QP$ such that kernel
  and cokernel are finitely generated projective modules over
  $C^*_rG\tensor \mathbb{K}$.)
\end{remark}

\section{Final remarks}

  We hope that Proposition \ref{Ktraceprop} can be applied to
  more situations than the one described in Theorem
  \ref{heredtrace}. However, in \cite{Schick(1999)} the
  situation where  $H$ is trivial (and consequently $X$ is finite) is
  classified. It turns
  out that in this setting the assumptions of Theorem
  \ref{Ktraceprop} can be fulfilled exactly if $G$ is a finite
  extension of a
  free group. But then one has a transfer homomorphism for the
  K-theory of the reduced $C^*$-algebras relating the trace for
  $G$ to the trace of the free subgroup of finite index. One easily
  computes the range of the trace using 
  this (and the known trace conjecture for the free group). The range
  is $\frac{1}{d}\integers$ where $d$ is
  the smallest index of a free subgroup. Therefore it is not
  necessary to give details
  of the  approach using Theorem \ref{Ktraceprop} which gives the
  same result.

\emph{Acknowledgments\/}. I am very much indebted to Warren
Dicks. Without his  help, I was
able to apply
Theorem \ref{Ktraceprop}  only in very basic situations.
Moreover, I thank Nigel Higson and John Roe
who pointed out that the
trace of a difference of two projections on a Hilbert space is an
index and therefore an integer and suggested that the same should work
in a more general setting, inspiring the proof of Lemma \ref{traceindex}.
I also thank  the referee for useful comments, in particular for pointing
out Remark \ref{remark:otherway}.


\end{document}